\documentclass[12pt]{article}
\usepackage{amssymb}

\usepackage{amsmath}
  \usepackage{paralist}
  \usepackage{graphicx}
  \usepackage{amssymb}
\usepackage[pagewise]{lineno}

\numberwithin{equation}{section}

\def \om{\omega}
\def \l{\lambda}
\def \beq{\begin{equation}}
\def \eeq{\end{equation}}
\def \la{\label}

\def \n{\noindent}

\def \a{\alpha}

\usepackage{color}

\def \R{\mathbb{R}}
\def \bs{\blacksquare}

\begin{document}
\title{On the Initial Value Problem for Hyperbolic Systems with Discontinuous Coefficients}

\author{Kayyunnapara Divya Joseph, \\
Department of Mathematics, \\
Indian Institute of Science Education and Research Pune, \\
Dr. Homi Bhabha Road, Pashan, \\
  Pune 411008, India.  \\
Email: divediv3@gmail.com} 

\maketitle

\numberwithin{equation}{section}
\numberwithin{equation}{section}
\newtheorem{theorem}{Theorem}[section]
\newtheorem{remark}[theorem]{Remark}

\begin{abstract}
Hyperbolic systems of the first and higher-order partial differential equations appear in many Multiphysics problems. We will be dealing with a wave propagation problem in a piece-wise homogeneous medium. Mathematically, the problem is reduced to analyzing two systems of partial differential equations posed on two domains with a common boundary. The differential equations may not be satisfied on the boundary (or part of the boundary), but some interface conditions are satisfied. These interface conditions depend on a specific physical problem. We aim to prove the existence and regularity of the solution for the case of hyperbolic systems of first-order equations with different domains separated by a hyperplane, where we need to formulate the interface conditions. We do this for the initial value problem in 1D-space variable when the coefficient matrix has discontinuity on $m$ lines.  More specifically, we find explicit solutions to the case when  the coefficient matrix is piecewise constant with a discontinuity along $1$ line or $2$ lines. We also prove the existence of solution to the general  initial value problem.  We then formulate the weak solution of  initial value problem for the corresponding symmetric hyperbolic system in $n $D-space variables  with interface conditions, get the energy estimates for this system, and prove the existence of solution to the system. 
\end{abstract}

\n { 2020 Mathematics Subject Classification. 35R05, 35L45, 35D30.}\\
\n {\bf Keywords:} multiphysics problem, hyperbolic system, discontinuous coefficient, interface conditions, existence and uniqueness theorems.

%\begin{small}

\section{Introduction} 

First-order linear and quasilinear hyperbolic equations appear in many fields of physics and engineering, such as continuum mechanics and relativity theory.  The most celebrated is the systems of compressible Euler's equations.  For the nonlinear case,  the existence of global-in-time solutions is well-understood only for simplified models in 1D - space variables. There is a well-developed theory for linear systems and higher-order scalar equations of hyperbolic type, with smooth coefficients. Rigorous study of such equations were started  by Hadamard \cite{h1}, continued by Schauder \cite{s1}, Petrowski \cite{p1}, Leray \cite{Leray}, Garding \cite{g1}, Friedrichs \cite{Fredrichs1}, \cite{Friedrichs2}, \cite{Friedrichs3}, Lax and Philips \cite{Lax2}, and others. Symmetric systems are most important because equations coming from mechanics and physics can be symmetrized according to \cite{Lax3}. 

In \cite{Hughes}, the authors establish well-posedness for a general class of quasi-linear evolution equations on a short time interval. They apply these results to second-order quasi-linear hyperbolic systems on $\R^n$ whose solutions $(u(t),u'(t))$ lie in the Sobolev space $H^{s+1}\times H^{s}$.  As a consequence, they obtain the well-posedness of the equations of elastodynamics if $s>2.5$, and of those of general relativity if $s>1.5$. 

The theory of initial boundary value problems for symmetric hyperbolic systems was developed by Friedrichs \cite{Friedrichs2}, Lax and Philips \cite{Lax2}, Rauch \cite{Rauch}, and for general hyperbolic systems, by Kreiss \cite{Kreiss} and Lopatinskii \cite{Ya}. They proved well-posedness and regularity in the Sobolev spaces $H^k =W^{k,2}$.

The problem considered in this paper is different from the classical theory described above, where a system with smooth coefficients is considered in a connected domain. Here, we consider two systems on two domains with a common boundary. The equation is not satisfied on the common boundary. Instead, some interface conditions are required which are given by physical applications. There is a large literature on this problem. See \cite{Christodoulou, Quarteroni} and reference therein, for different types of equations coming from applications. In our paper, we are not dealing with numerical methods that are well-developed but only analytical methods. For these numerical methods, see \cite{Fatemeh} and the references therein, though \cite{Fatemeh} is also dealing with some theoretical aspects. See also NASA report \cite{Kopriva} for numerical methods, along with some theoretical issues. There, a multidomain Chebyshev spectral collocation method is developed to solve hyperbolic partial differential equations. The results showing the performance of this method on one-dimensional linear models, and one- and two-dimensional nonlinear gas-dynamics problems are presented.

In this work, we have rigorously proved the existence of a Lipschitz continuous solution of the 1-dimensional system which is $C^1$ away from the lines of discontinuities of the coefficient matrix. We also considered the n-dimensional system and proved the existence of a unique weak solution. We derived the interface condition for the system. Such an analysis is not done previously as per our knowledge. \\

We have a two-fold motivation for our study. As a multiphysics problem, we chose acoustic boundary-contact problems \cite{Boris}.  By multiphysics problems, we mean problems analyzing multiple, simultaneous physical phenomena.  As an example of the methodology, we choose the papers \cite{Friedrichs2, Lax2, Rauch, Fatemeh}. The authors of \cite{Boris} are dealing with a multiphysics problem but this old study was prepared before the multiphysics methodology was developed. Using modern language, the model studied in \cite{Boris} describes two systems on two domains with a common boundary.  The acoustic medium is assumed to be vertically stratified and bounded by a plate with concentrated inhomogeneities. Mathematically, the problem is formulated as the Helmholtz equation with discontinuous coefficients, which is \eqref{1.2} below. It originates from the wave equation
\begin{equation}
v_{tt} =  c^2(x,z) ( v_{xx} + v_{zz} )
\label{1.1}
\end{equation}
where $c(x,z)$ is positive. The physical meaning of $c(x,z)$ is the speed of acoustic waves and it is assumed in \cite{Boris} that the acoustic medium is layered so that $c$ is a function of $z$ alone that has a discontinuity along $z = H$. This equation appears in many multiphysics problems \cite{Boris, Fatemeh}. Without loss of generality, we can assume $H=0$. Taking $v(x, z, t)=e^{\pm i \om t}g(x, z), \om \neq 0$ we get 
\begin{equation}
\Big( g_{xx} + g_{zz} +  \frac{\om^2}{c^2(z)}\Big) g(x, z)=0,
\label{1.2}
\end{equation}
This equation is called the Helmholtz equation and is typical for the problems of wave propagation. 

More recently, Fatemeh et al. \cite{Fatemeh} study the system of equations in 1D - space variable 
\begin{equation}
\begin{aligned}
&u_t + A u_x = 0, -1 < x < 0, t > 0, \\
 &u(x, 0) = f(x), -1 < x < 0, t = 0,
\label{1.3}
\end{aligned}
\end{equation}
and 
\begin{equation}
\begin{aligned}
 &v_t + B v_x = 0, 0 < x < 1, t > 0, \\
 &v(x, 0) = g(x), 0 < x < 1, t = 0,
\label{1.4}
\end{aligned}
\end{equation}
where  $A, B$ are $m \times m$ and $n \times  n$ symmetric constant matrices, respectively, and $u, v$ are unknown vectors of sizes $m$ and $n$ respectively. In general $m \neq n$. They ignore the boundary conditions at $x= \pm 1$ and focus only on the interface conditions at $x=0$.   Their work was on the numerical approximations of the solution to the initial value problem to \eqref{1.3}, \eqref{1.4}. They study the interface conditions for the problem to be well-posed as well as the stability conditions at the interface. They prove the stability of a semidiscrete approximation of the solution. \\

The second order wave equation \eqref{1.1} with variable coefficients may be converted to a symmetric hyperbolic system of the form,  
\begin{equation}
B_0(z) u_t + B_1(z) u_x + B_2(z) u_z = 0
\label{1.5}
\end{equation}
where 
\begin{equation*}
\begin{aligned}
&B_0(z) = 
\begin{bmatrix}
c^2(z) & 0 & 0\\
0 & c^2(z) & 0\\
0 & 0 & 1
\end{bmatrix},
B_1(z)  = 
\begin{bmatrix}
0 & 0 & - c^2(z) \\
0 & 0 & 0\\
- c^2(z) & 0 & 0
\end{bmatrix}, 
\end{aligned}
\end{equation*}
\begin{equation*}
\begin{aligned}
B_2(z)  = 
\begin{bmatrix}
0 & 0 & 0\\
0 & 0 & - c^2(z) \\
0 & - c^2(z)  & 0 .
\end{bmatrix}
\end{aligned}
\end{equation*}
Here $u = (v_x, v_z, v_t)^T$ is a 3D vector, with initial condition, $u(x, z, 0) = u_0(x, z).$ Here the upper index $^T$ denotes transposition.

The reason for the choice of these functions as unknown variables follows from the weak formulation of
the wave equation (\ref{1.1}) %$v_{tt}=v_{xx} +v_{zz}$ 
$$\int_{\Omega} v_t V_t d\Omega=\int_{\Omega} c^2(z)(v_x V_x+v_z V_z)d\Omega$$
where $V$ is a smooth function supported in the space-time region $\Omega$, so that, $V=0$ on the boundary $\partial\Omega$.

Since ${c^2}\ne 0$ the system (1.5) is equivalent to  
\begin{equation}
u_t + \tilde B_1(z) u_x + \tilde B_2(z) u_z = 0
\label{1.6}
\end{equation}
where
\begin{equation*}
\tilde B_1(z) = 
\begin{bmatrix}
0 & 0 & -1\\
0 & 0 & 0\\
-c^2(z) & 0 & 0
\end{bmatrix}
\text{and} \,\,\,  \tilde B_2(z) =
\begin{bmatrix}
0 & 0 & 0\\
0 & 0 & -1\\
0 & -c^2(z) & 0.
\end{bmatrix}
\end{equation*}

{ \flushleft \bf Definition: } (see L. C. Evans \cite{Evans}) \\
(i) We say that the system 
\begin{equation}
u_t + \tilde B_1(x, z, t) u_x + \tilde B_2(x, z, t) u_z = 0
\label{1.7}
\end{equation}
 is symmetric hyperbolic if the matrices $ \tilde B_j(x, z, t),\;j=1,2$ are symmetric $\forall x, z \in \R, \,\,\,  t \geq 0.$    \\
(ii) We say that this system is  strictly hyperbolic if $\forall x, z \in \R,\;{\eta} \in \R^2$ with $|\eta| =1$ and $\forall t \geq 0$ the matrix $B(z, t; \eta)$ given by
\begin{equation*}
B(z, t; \eta) = \sum_{j=1, 2} {\eta}_j \tilde B_j(x, z, t),  \,\,\, x, z \in \R, \,\,\, t \geq 0.
\end{equation*}
 has distinct real eigenvalues. The symmetry of $B_1$ and $B_2$ implies that the matrix $B(z, t; \eta)$ has real eigenvalues and a complete set of eigenvectors for each $|\eta|=1$. This implies that the system is hyperbolic. Hence, symmetry of $\tilde{B_1}(x,z,t)$ and $\tilde{B_2}(x,z,t)$ in \eqref{1.7} implies hyperbolicity, see \cite{Evans}. 

Before we proceed further, we note that the function $c(z)$ is assumed to be positive $\forall z$.
{ \flushleft \bf Lemma 1.1: } 
Let $c^2(z)\not\equiv 1$. Then: 

(i) The system \eqref{1.6} is not symmetric hyperbolic.

(ii) The system \eqref{1.6} is strictly hyperbolic.  \\
%\end{lemma}
{  \bf Proof:} (i) It is easy to see that this system is not symmetric hyperbolic unless $c^2(z)\equiv 1$ in which case the matrices $\tilde B_j,\,j=1,2$ are symmetric.

(ii) By direct calculation, we find that the eigenvalues of the matrix 
\[
\eta_1 \tilde B_1(x,z,t)+ \eta_2 \tilde B_2(x,z,t), \,\,\,\text{with} \,\,\, \eta_1^2+\eta_2^2 =1
\]
are  $- c(z), 0,  c(z)$. Since $c(z)>0$, they are real and distinct $\forall z$, and the system \eqref{1.6} is strictly hyperbolic.   $\bs$ 

\begin{remark}We note that  the matrices $\tilde{B_1}, \tilde{B_2}$ have a discontinuity along $z=0$. %\\
\end{remark} 

We now describe a plan for the remaining part of the paper. First, we study \eqref{1.7} in the 1D - space. In this case, the system becomes
\begin{equation}
u_t + B_1(z, t) u_z = 0 .
\label{1.8}
\end{equation}

In Section 2, we consider a more general system than \eqref{1.3}-\eqref{1.4} or  \eqref{1.5} of $n$ equations in $u$ having $n$ components, with the corresponding matrix having a discontinuity along $z=0$. We prove the existence and uniqueness of the solution to the Cauchy problem for this system. In Section 2, we first consider the initial value problem for the system in 1D - space variable in the domain $\{(z,t): z \in \R^1, t>0\}$ with $B(z,t)$ having discontinuities along $m$ lines $z=z_l, \,\,\, l = 1, \dots, m$ and study the existence of the solution with the interface condition along these lines. In Section 2.1 we construct explicitly the solution to the initial value problem for \eqref{1.8} when $B$ is a piecewise constant matrix with a discontinuity only along $\{(z,t): z=0\}$ using the method of characteristics. We also do the same for $B$ having  2 discontinuities along $z=z_1, \,\,\, z=z_2$. In Section 2.2 we prove the existence of a unique solution, to the initial value problem to the general system \eqref{1.8}, which is once continuously differentiable except along the line of discontinuity of $B$, with an interface condition on $\{(z,t): z=z_l,\;t>0, \,\,\, l = 1, \dots, m\}$. We follow Friedrichs \cite{Friedrichs3} (see also Courant and Hilbert \cite{Courant} and  John \cite{John}.) In all these works,  the method of characteristics is used for systems having smooth coefficients, but in our work, we use this method of characteristics for a system having discontinuous coefficients, with a different kind of analysis, and achieve different results.   The interface condition for the solution constructed in Section 2.2 is also derived here, explicitly. 

In Section 3, we derive energy estimates that we use to construct solutions to the initial value problem for the $n$ - space dimensional symmetric hyperbolic system with discontinuity along the surface $\{(x_1,\dots,x_n): x_n=0\}$ and derive the corresponding interface conditions.

We would like to clarify the connection between sections 2 and 3. If the matrix $B$ in section 2 is symmetric, we can apply the same energy method that we used in section 3, but we would get a weaker result. The solution obtained by the characteristic method is in the space of continuously differential functions except along the lines of discontinuity of the coefficient matrices but solutions obtained by the energy method are less regular and only weakly differentiable.  Further if the number of the space variables, $n\ge 2$, we cannot use the method of characteristics.

\section{One space variable case: method of characteristics}
\la{St-L_sect}

In this section, we consider the initial value problem for the system in 1D - space variable in the domain $\{(z,t): z \in \R^1, t>0\}$ 
\begin{equation}
\begin{aligned}
&u_t  + B(z,t) u_z = 0,\,\,\,   z \in (-\infty, z_1) \cup ( \cup_{i=1}^{m-1} (z_i, z_{i+1} ) ) \cup (z_m,\infty), t>0\\ 
&u(z, 0) = u_0(z),\,\,\, z \in \R^1
\end{aligned}
\label{2.1}
\end{equation}
Our focus is on the existence of the solution with the interface condition on the lines $z=z_l, \,\,\, l = 1,\dots, m$, where $B(z,t)$ has a discontinuity. For a 1D - problem, the solution space can be chosen to be a classical space of continuously differentiable functions with $C^1$ norm. It is well known that even for the 1D wave equation, the solution has the same regularity as the initial data in classical derivatives, i.e the solution is $C^1$ if the initial data is $C^1$. The method of characteristics, which typically works for systems in 1D - space variable, gives regularity in the classical sense since the solution space can be chosen as $C^1$ space, see \cite{John, Friedrichs3}. Yet, though there is a loss of regularity for the higher dimensional wave equation, regularity in Sobolev spaces is preserved, see John \cite{John} and Evans \cite{Evans}. So, we use the method of characteristics rather than the energy method. We use the energy method later for the 2D - space variable system \eqref{1.7}. 

First, we recall definitions of characteristic speed, characteristic direction, and characteristic curve which are central to our analysis. We always assume the system is strictly hyperbolic.

{\flushleft \bf Definition:} 
(i) The eigenvalue $\lambda_j, j=1,\dots, n$ of  $B$ is called the $j^{\rm th}$ - characteristic speed for the system \eqref{2.1}.
 
{\flushleft 
(ii)} The curve with speeds $\lambda_j$ passing through $(z,t)$, defined by
\[
\frac{d \alpha}{ds} =\lambda_j(\alpha,s),\,\,\, \alpha(s,z,t)|_{s=t}=z
\]
is called the $j^{\rm th}$ - characteristic curve through  $(z,t)$. 

 The importance of $j^{\rm th}$ - characteristic speed is that it is the speed of propagation of the component of the initial data in the direction of $j^{\rm th}$ eigenvector of $B$. We need the following smoothness properties in future analysis. We consider a more general hyperbolic system \eqref{2.1} than the examples above in Section 1. We introduce the following assumptions.
{\bf  \flushleft Assumptions on $B(z, t)$: } \\
H1: Assume that $B(z,t)=(b_{ij}(z,t))$ is an $n \times n$  matrix with $b_{ij}$ and its derivatives up to order two are bounded except along $m$ lines $z=z_l, l=1,\dots,m$ in the $(z,t)-$plane and that the traces $(\partial_t)^l(\partial_z)^m b_{ij} (z_l\pm,t)$ exist for $l+m\leq 2$.  \\
H2:  $B(z,t)$ has real distinct eigenvalues,
\[
{\l}_{1}(z,t)< {\l}_{2}(z,t)< \dots <  {\l}_{k}(z,t)<0
< {\l}_{k+1}(z,t)< \dots < {\l}_n(z,t).
\]

Since the eigenvalues are distinct, $B$ is diagonalizable \cite{Dym}, \cite{Hoffman}. The next lemma is on the properties of the diagonalizing matrix.
 
{\flushleft \bf Lemma 2.1:} Let $B$ be $B(z, t)$ from the system \eqref{2.1}. There exists an $n\times n$ matrix $A=A(z,t)$  with the following properties.  
\\
(i) $A$ and its inverse $A^{-1}$ are twice continuously differentiable except on  $z= z_l$  \\
(ii) The derivatives of $A, A^{-1}$ up to order 2 are bounded on every bounded set in $z \neq z_l$ with right and left traces  at $z=z_l$, for $l=1,\dots,m$, and   
\begin{equation}
B = A \Lambda A^{-1} ,\,\,\,  \Lambda = diag( \l_1,\dots,\l_n),
\label{2.2}
\end{equation}
 $\l_j$ are $C^2$ for $z \neq z_l$ and $(\partial_t)^i(\partial_z)^j {\l}_j (z_l\pm,t)$ exists. In fact $A = (\bar{r}_1,\bar{r}_2,..., \bar{r}_n)$ where  $\bar{r}_j$ is eigenvector corresponding to speed $\l_j$. 
\\
{\bf Proof :}
According to $H_2$, $B(z, t)$ in \eqref{2.2} has distinct eigenvalues and hence simple roots of the characteristic polynomial, whose coefficients are polynomials in $b_{ij}$ are clearly twice continuously differentiable on $(z,t)$ except $z=z_l$. This gives $\l_j$ are $C^2$ and, as in Lax \cite{Lax},  the corresponding eigenvectors may be chosen in $C^2$. Now, $A$ is formed by these eigenvectors as column vectors. So, $A$ is $C^2$ except $z=z_l$. Now, to show that $\lambda_j$ is uniformly bounded, we take an eigenvalue $\l$ of $B,$ with corresponding eigenvector $\bar{r}=(r_1,r_2,...r_n)$. Note that  
$\l <\bar{r}, \bar{r}> = <B\bar{r}, \bar{r}>$, which gives
\begin{equation*}
\begin{aligned}
|\l | \leq \frac{ | <B\bar{r}, \bar{r}> | }{| <\bar{r}, \bar{r}> |} = \frac{ | \sum_{i, j = 1}^{n}  b_{ij} \,\,\, r_i \,\,\, r_j |}{  \sum_{i, j = 1}^{n}  {| r_i| }^2  }
&\leq  \sum_{i, j = 1}^{n} | b_{ij} | \frac{| r_i | }{ ||\bar{r}||} \frac{| r_{j} | }{ ||\bar{r_j}||}\\
&\leq  \sum_{i, j = 1}^{n} | b_{ij} | \\
&\leq n^2\sup_{\R \times [0,\infty)}{ | b_{ij}(z, t) | }. 
\end{aligned}
\end{equation*}
This shows $\sup_{\R\times [0,\infty)} |\lambda(z,t)| \leq n^2 \sup_{\R \times [0,\infty)}{ | b_{ij}(z, t) | }$, which is bounded
by $H_1$. 

Also, the structure of the inverse matrix $A^{-1}$ gives that (a) $A^{-1}$ is also $C^2$ in $z \neq z_l$; (b) its derivatives of order up to 2 are bounded on  $\R\times [0, T]$, $\forall \,\,\, T>0$. Now, our assumption that the right and left traces of $B$ at $z=0$ exist implies that the right and left traces of $A, A^{-1}$ and their derivatives of order 2 at $z=0$ also exist. $\bs$

\subsection{Piecewise constant case}
First, we take the simplest case when the matrix $B$ is piecewise constant with a discontinuity along $z=0$. Thus, we assume
\begin{equation} 
B(z) = \begin{cases}
B_{+} &\text{ if } z >0 \\
B_{-} &\text{ if } z<0
\end{cases}
\label{2.3}
\end{equation}
where $B_{+}, B_{-}$  are constant $n \times n$ matrices, $B_+\ne B_-$. Physically, that means that we consider the wave propagation in two different homogeneous media occupying the half-spaces $z<0$ and $z>0$.
 
Introduce the continuity interface condition, 
\begin{equation}
(A^{-1}u)(0+, t)=(A^{-1}u)(0-, t). 
\label{2.4}
\end{equation}

{\bf  \flushleft  Theorem 2.2:} { Assume that the initial data $u_0(x)$ is continuously differentiable. Then under} the  assumptions \eqref{2.3} on $B$, there exists a unique continuously differentiable solution $u(z,t)$ to \eqref{2.1} for $z \neq 0$. Further, it satisfies the interface condition \eqref{2.4}, i.e., $A^{-1}u$ is continuous across $z=0$.  \\

{ \bf Proof :} Let $u = A v$ for $z \neq 0.$ Then the system of partial differential equations in \eqref{2.1} becomes 
$v_t + A^{-1} B(z) A v_z =0,$ which in turn gives 
\[
v_t + \Lambda v_z =0,\,\,\, 
\] 
where $\Lambda= diag(\lambda_1,\dots,\lambda_n)$, for $z \neq 0.$ Let $v=(v_1, \dots, v_n)^T$, and $A^{-1} u_0(z)=(v_{0 1}(z),\dots,v_{0 n}(z))^T $. Then the problem for $u$ reduces to  
\begin{equation}
(v_j)_t  + \l_j (v_j)_z = 0, \,\,\,
(v_j)(z, 0) = v_{0 j}(z) .
\label{2.5}
\end{equation}

Since $B(z)$ is constant on $z>0$ and $z<0$, eigenvalues and eigenvectors also have the same property. In the following, let us denote  by $\lambda^{+}$ the value of $\lambda(z)$ for $z>0$ and  by  $\lambda^{-}$ the value of $\lambda(z)$ for $z<0$. We construct the solution in a different way for characteristic speed $\l_j>0$ and $\l_j<0$.
\\
{\bf The case $\l_j >0 :$ } 
\begin{figure}[!hbtp]
\includegraphics[width=13cm,height=4cm]{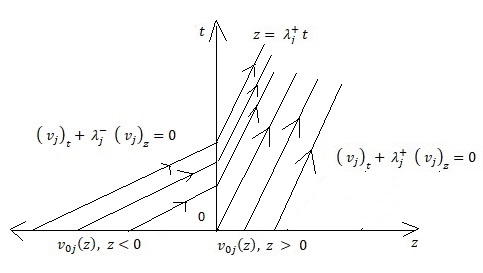}
\caption{Characteristics for $\l_j>0$}
\end{figure} 
We first solve \eqref{2.5} in the region $z < 0.$ 
Using the method of characteristics we get the characteristic curve 
$z = \l_j^{-} t +z_0,$ where $z_0$ is an arbitrary constant and  the solution to \eqref{2.5} is
\[
v_j(z, t) =  v_{0 j}(z_0)= v_{0 j}(z - \l_j^{-} t), z<0.
\]
Here we use the fact that the solution is constant along the characteristic curves. To find the solution  for $z>0$ observe that $v_j(0 -, t)$ exists. If $z \geq \l_j^{+} t$ the solution is given by the method of characteristics and data coming from $t = 0,$ namely
\[
v_j(z, t) = v_{0 j}(z - \l_j^{+} t), z \geq \l_j^{+} t.
\]
To find the formula in the region $0< z<  \l_j^{+} t$
draw the characteristic from $(x, t)$ backwards in time to hit $z=0$ at $(0,  \tau)$ for some $\tau$. The characteristic is given by $ \frac{z - 0}{t - \tau} = \l_j^{+},$ which results in $\tau = t - \frac{z}{\l_j^{+}}.$ so that
\[
v_j(z, t) = v_j(0 - , \tau) = v_{0 j}( - \l_j^{-} \tau ) = v_{0 j}\Big( \frac{ \l_j^{-}}{\l_j^{+} } (z -  \l_j^{+} t ) \Big),
\] 

for $0< z<  \l_j^{+} t.$  So, the solution  for speed  $\lambda_j>0$, is   
 \begin{equation}
v_j(z, t) = \begin{cases}
v_{0j} (z - \l_j^{-} t ) &\text{if } z \leq 0, \\
v_{0j}( \frac{\l_j^{-}}{\l_j^{+}} (z -  \l_j^{+} t)), &\text{if } 0< z< \l_j^{+} t, \\
v_{0j} (z - \l_j^{+} t), &\text{if } z \geq \l_j^{+} t.
\end{cases}
\label{2.6}
\end{equation}

{\bf The case $\l_j <0 :$} First, we consider the region $z > 0$.
The characteristic passing through $(z, t)$ is of the form $z = \l_j^{+} t +z_0.$ The equation $(v_j)_t  + \l_j^{+} (v_j)_z = 0$
states that the solution is constant along characteristics.
So, the solution  is $v_j(z, t) =  v_{0 j}(z_0)= v_{0 j}(z - \l_j^{+} t), z>0.$ 
Hence, $v_j(0 +, t)$ exists. %We then solve in the region  $z<0.$ 
\begin{figure}[!hbtp]
\includegraphics[width=13cm,height=4cm]{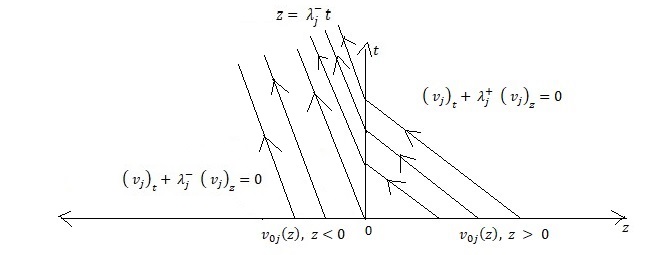}
\caption{Characteristics for $\l_j<0$} 
\end{figure}
We solve the problem $(v_j)_t  + \l_j^{-} (v_j)_z = 0, \,\,\, (v_j)(z, 0) = v_{0 j}(z), z< 0.$ In the region
$z \leq \l_j^{-} t$, the solution is given by the method of characteristics and data coming from $t = 0,$ 
\[
v_j(z, t) = v_{0 j}(z - \l_j^{-} t), z \leq \l_j^{-} t.
\]
In the region $  \l_j^{-} t<z<0$,  the speed is $ \l_j^{-}.$ We draw the characteristic from $(x, t)$ backward in time to hit $z=0$ at $(0,  \tau)$, for some $\tau$.
The value $(v_j)(0, \tau)$  at that point is determined by solving 
$(v_j)_t  + \l_j^{+} (v_j)_z = 0$. Note that the characteristic is given by $ \frac{z - 0}{t - \tau} = \l_j^{-},$ which gives $\tau = t - \frac{z}{\l_j^{-}}.$ So, the solution is 
\[
v_j(z, t) = v_j(0 + , \tau) = v_{0 j}( - \l_j^{+} \tau ) 
= v_{0 j}( \frac{ \l_j^{+}}{\l_j^{-} } \Big(z -  \l_j^{-} t ) \Big),  \l_j^{+} t < z< 0.
\]
Combining the results for all the regions, we find the solution for speed  $\lambda_j<0$ 
\begin{equation}
v_j(z, t) = \begin{cases}
v_{0j} (z - \l_j^{-} t ) &\text{if } z \leq \l_j^{-} t, \\
v_{0j}( \frac{\l_j^{+}}{\l_j^{-}} (z -  \l_j^{-} t)), &\text{if } \l_j^{-} t < z< 0, \\
v_{0j} (z - \l_j^{+} t), &\text{if } z \geq 0.
\end{cases}
\label{2.7}
\end{equation}
Hence, the solution for \eqref{2.1} is $u=A v,\;{\rm where}\;v=(v_1,\dots,v_n)$ and $v_j,j=1,\dots,n$ is given by \eqref{2.6} and \eqref{2.7}. Since $v$ is continuous across $z=0$ and $v=A^{-1} u,$ we conclude that $A^{-1} u$ is continuous across $z=0$, and so the interface condition $A^{-1} u(0-, t) = A^{-1} u(0+, t)$ is satisfied. \\
 $\bs$

The same method as above can be used to construct a solution to the initial value problem with discontinuities on $m>1$ lines: $z=z_l, l=1,\dots,m$. We write down the solution for the case $m=2$. Given an eigenvalue, it must have the same sign in every region by the assumption H2. Otherwise, the solution will be in the space of measures, see \cite{KTJ}, \cite{Lefloch}. We will not be considering eigenvalues having different signs in this case. \\

{\bf Solution when the coefficient matrix $B$ is piecewise constant with two lines of discontinuities, i.e.}
\begin{equation} 
B(z) = \begin{cases}
B_{1} &\text{ if } z <z_1 \\
B_{2} &\text{ if } z_1<z<z_2 \\
B_{3} &\text{ if } z>z_2
\end{cases}
\label{2.8}
\end{equation}
 where $B_{i}, i=1,2,3$  are constant $n \times n$ matrices. Let us denote the eigenvalues of the matrix $B_1,B_2,B_3$ in \eqref{2.8} by $\l_{1j}, \l_{2j}, \l_{3j}, j =1,\dots, n$ respectively. Following the same analysis as above, we can get the following formula for the component in the $j^{\rm th}$ characteristic direction $v_j,  j =1, 2, \dots, n$ in the following form. The solution $u(z, t)$ to system \eqref{2.1} for $z \neq z_1, z \neq z_2$ satisfying interface condition $A^{-1}u({z_i}+)=A^{-1}u({z_i}-), i =1, 2$ is $u = Av, v = (v_1,\dots, v_n) $ with $v_j$ given by \eqref{2.9} below when the eigenvalues are all positive and by \eqref{2.10} below when the eigenvalues are all negative.  \\
Case  1: Let $\l_{1j}, \l_{2j}, \l_{3j}$ be all positive: (Refer to Figure 3.)
\begin{figure}[!hbtp]
\includegraphics[width=15cm,height=5cm]{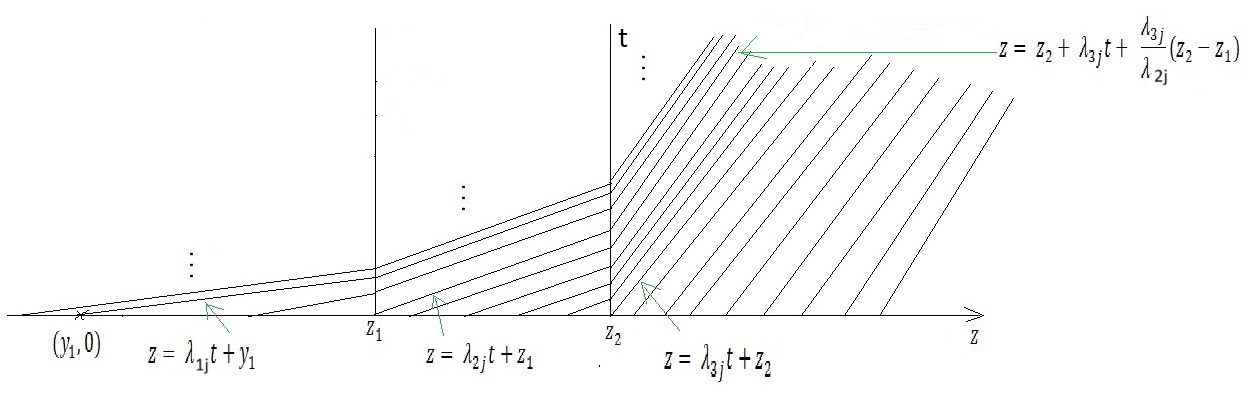}
\caption{$\l_{1j}, \l_{2j}, \l_{3j}$ are all positive}
\end{figure} 
\begin{equation}
v_j(z, t) = \begin{cases}
v_{0j} (z - \l_{1j} t ) &\text{if } z \leq z_1, \\
v_{0j} (z - \l_{2j} t ) &\text{if }  z_1\leq z \leq z_2, z\geq \l_{2j} t +z_1 \\
v_{0j}( z_1+\frac{\l_{1j}}{\l_{2j}} (z -  \l_{2j} t) -z_1), &\text{if }, z_1<z<z_2,  z< \l_{2j} t +z_1, \\
v_{0j} (z - \l_{3j} t), &\text{if } z \geq \l_{3j} t +z_2,\\
v_{0j}( z_1+\frac{\l_{1j}}{\l_{3j}} (z -  \l_{3j} t) \\
+\frac{\l_{1j}}{\l2j}(z_2-z_1)), &\text{if } z\geq z_2, z \leq z_2+\l_{3j}t  +\frac{\l_{3j}}{\l_{2j}} (z_2-z_1), \\
v_{0j}( z_2+\frac{\l_{2j}}{\l_{3j}} (z -  \l_{3j} t) -z_2), &\text{if } z\geq z_2, z_2+\l_{3j}t  + \frac{\l_{3j}}{\l_{2j}} (z_2-z_1)\\
&\leq z< \l_{2j} t +z_1. \\
\end{cases}
\label{2.9}
\end{equation}
 for j =1, 2, \dots, n. \\
Case 2: Let $\l_{1j}, \l_{2j}, \l_{3j}$ be all negative: (Refer to Figure 4.)
\begin{figure}[!hbtp]
\includegraphics[width=15cm,height=5cm]{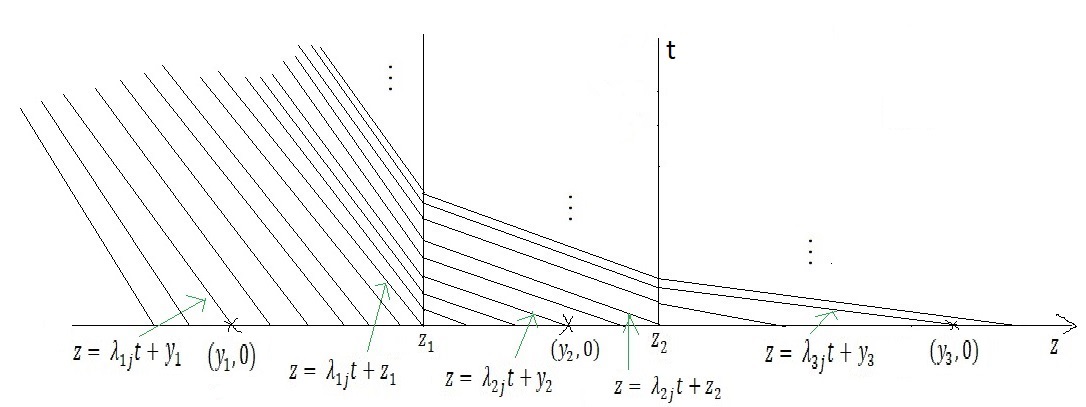}
\caption{$\l_{1j}, \l_{2j}, \l_{3j}$ are all negative.}
\end{figure} 
\begin{equation}
v_j(z, t) = \begin{cases}
v_{0j} (z - \l_{3j} t ) &\text{if } z \geq z_2, \\
v_{0j} (z - \l_{2j} t ) &\text{if }  z_1\leq z \leq z_2, z\leq \l_{2j} t +z_2 \\
v_{0j}( z_2+\frac{\l_{3j}}{\l_{2j}} (z -  \l_{2j} t) -z_2), &\text{if }, z_1<z<z_2,  z \geq \l_{2j} t +z_2, \\
v_{0j} (z - \l_{1j} t), &\text{if } z \leq \l_{1j} t +z_1,\\
v_{0j}( z_2+\frac{\l_{3j}}{\l_{1j}} (z -  \l_{1j} t-z_1)\\
 +\frac{\l_{3j}}{\l2j}(z_2-z_1)), &\text{if } z_1 \geq z \geq z_1+\l_{1j}t+\frac{\l_{1j}}{\l_{2j}}(z_1-z_2), \\
v_{0j}( z_1+\frac{\l_{2j}}{\l_{1j}} (z -  \l_{1j} t) -z_1), &\text{if } z\leq z_1, \l_{1j} t +z_1 < z\\
& \hspace{0.3cm}<   z_1+\l_{1j}t  +\frac{\l_{1j}}{\l_{2j}}(z_2-z_1). \\
\end{cases}
\label{2.10}
\end{equation}
 for j =1, 2, \dots, n. \\
 $\bs$
\subsection{Solution for the general case  of system \eqref{2.1}} 
In this section, we consider the general initial value problem \eqref{2.1} with $B(z,t)$ satisfying hypothesis $H_1$ and $H_2$. We shall prove the following theorem. \vspace{0.5 cm}

{\bf  \flushleft Theorem 2.3.} Let the initial data $u_0$ be continuously differentiable. Under the assumptions $H_1$ and $H_2$ on $B$ there exists a unique solution $u(z,t)$ to \eqref{2.1} in the space of continuously differentiable functions except at  $z = z_l, l=1,\dots,m$ with interface condition \eqref{2.4}. 
{   \flushleft First, we need a few lemmas. }
{\bf  \flushleft Lemma 2.4:} 
Let $\l(z, t)$ be a bounded continuous function with bounded first order derivatives in $z \neq z_l, \,\,\, l=1,\dots,m$ with $z_1 < z_2 < \dots < z_m$. Assume 
$(\partial_z)^i\l(z_l+, t), \,\,\,(\partial_z)^i \l(z_l-, t),$ $i=0,1$ exists, then,\\ 
(i) The initial value problem
\begin{equation}
\begin{aligned}
&\frac{d \a}{d s} = \l(\a(s, z, t), s), \,\,\,  0<s<t \\
& \a(s=t, z, t)= z
\end{aligned}
\label{2.11}
\end{equation}
has a Lipschitz continuous solution $ \a(s, z, t),$ which is continuously differentiable in variables $s, z, t$ except on the line $z=z_l,$ and its first order derivatives have a left and right limit along $z=z_l$. \\
(ii) Further, let $\tau_l(z,t)$ be the point where $\alpha(\tau_l(z,t),z,t)=z_l$. Then $\tau_l$ is a continuously differentiable function of $(z,t)$ for $z \neq z_l$. In fact, $ \a(s, z, t)$ is twice  continuously differentiable in $s,$ except when $s \in \{ \tau_1, \tau_2, \dots \tau_m  \}.$  \\
{\bf Proof: } 
 By the Fundamental Theorem of Calculus and classical fixed point argument, the problem \eqref{2.11} is equivalent to 
\begin{equation}
\alpha(s,z,t)=z + \int_{t}^s \lambda(\alpha(s_1,z,t), s_1) d{s_1}
\label{2.12}
\end{equation}
and has solution  $ \a(s, z, t)$ which is Lipschitz continuous in all the variables. Indeed, the Lipschitz continuity in $s$ follows from
\begin{equation*}
\begin{aligned}
|  \a(s_2, z, t) -  \a(s_1, z, t)  | &= | \int_{s_1}^{s_2}{  \l( \a(s, z, t), s)  } ds |  \\
&\leq ||\lambda||_\infty | \int_{s_1}^{s_2}{dt} |
=  ||\lambda||_\infty |  s_2 - s_1 |
\end{aligned}
\end{equation*}
So, $\a(s, z, t)$ is Lipschitz continuous in $s$. 

[ In this part we use  Hartman  (Theorem 3.1, page 95) \cite{Hartman}. $y_0, t_0$ in this theorem is $z, t$ respectively and $y$ in this theorem is $\a.$ I have also included more details here. ] \\
 
If $\lambda<0$, then $\alpha(s,z,t)$ is a decreasing function of $s$, for all $0\leq s \leq t$. So if $z>z_m$ then $\alpha(s,z,t) > z_m$, for all $0<s<t$ and $\alpha(s=t,z,t)=z$. If $\a(s,z,t)$ does not touch the line $z=z_m$ like in Figure 5 below, then we get $\a$ is $C^1$ in the variable $s$ for $0 \leq s \leq  t$. We note that $(z, t)$ is the initial data which we consider to be in the region $z>z_m$ in Figure 5 and the region $z_{m-1}< z< z_m$ in Figure 6. $z$ is seen on the horizontal axis whereas the parameter $s$ is on the vertical axis, varying on $[0, t].$  The curve in Figure 5 and Figure 6 is a graph of $\alpha(s,z,t)$ as a function of the parameter $s$. 

\begin{figure}[!hbtp]
\includegraphics[width=15cm,height=5cm]{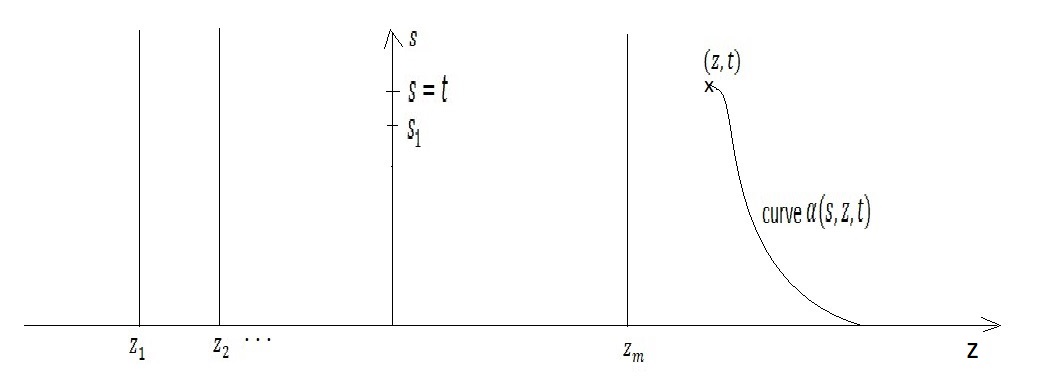}
\caption{ When  $\a(s,z,t)$ does not touch the line $z=z_m$  }
\end{figure} 

If $\a(s,z,t)$ cuts the line $z=z_m$ at ${\tau}_m(z, t)$ then \eqref{2.12} gives us that $\a$ is $C^1$ if ${\tau}_m(z, t)< s \leq t,$ like we see in the Figure 6.\\
\begin{figure}[!hbtp]
\includegraphics[width=15cm,height=5cm]{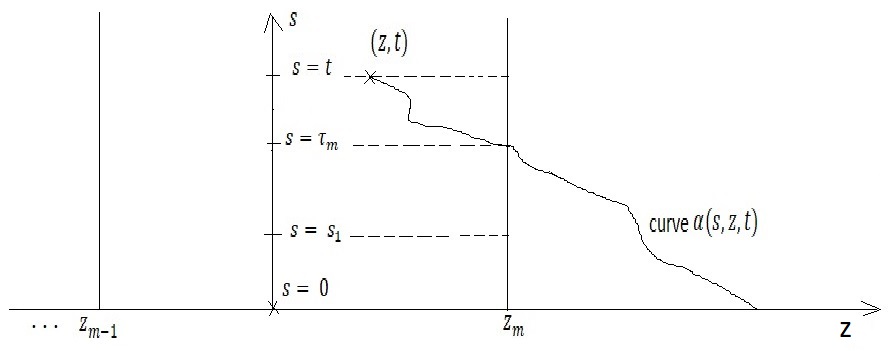}
\caption{ $\a(s,z,t)$ cuts the line $z=z_m$ at ${\tau}_m(z, t)$ }
\end{figure} 

We now show that $\a$  is $C^1$ in $s$ on the remaining curve $0 < s < \tau_m(z, t).$ For $s< \tau_m(z, t),$ we have
\begin{equation}
\begin{aligned}
\a(s,z,t) &= \int_t^s{ \l(\a(s_1,z,t), s_1) d{s_1}} + z \\
&= z +  \int_t^{\tau_m(z, t) }{ \l(\a(s_1,z,t), s_1) d{s_1}} + \int_{\tau_m(z, t) }^s{ \l(\a(s_1,z,t), s_1) d{s_1}}\\
&=c(z, t) +  \int_{\tau_m(z, t) }^s{ \l(\a(s,z,t), s_1) d{s_1}}
\label{2.13}
\end{aligned}
\end{equation}
where $c(z, t)=z +  \int_t^{\tau_m(z, t) }{ \l(\a(s,z,t), s_1) d{s_1}}$ is clearly constant in the variable $s.$  We can see in  Figure 6  that, for $s_1 <\tau_m(z, t), \,\,\, \a(s_1,z,t) > z_m,$ so  $ \l(\a(s_1,z,t), s_1)$ is continuous in $s_1,$  for $0< s_1 <\tau_m(z, t)$ since $\a$ is Lipschitz continuous in $s_1.$  Using  \eqref{2.13} $\a$  is $C^1$ in $s$ in the interval $0< s <\tau_m(z, t).$ Differentiating \eqref{2.13} we have $\frac{d \a}{d s} = \l(\a(s, z, t), s)$ which is $C^{1}$. So   $\a$ is  $C^{2}$ in the variable $s,$  in the  interval $0< s< {\tau}_m(z, t).$ By a similar argument,  $\a$ is  $C^{2}$ in the variable $s,$  in the other interval $\tau_m(z, t)< s \leq t.$ So we have that $\a$ is  $C^{2}$ in the variable $s,$ except on $s = \tau_m(z, t),$ when initial data $(z, t)$ is such that $z_{m-1}< z< z_m.$ We can do the same when the initial data $(z, t)$ is such that $z_{l-1}< z< z_l, \,\,\, l = 1, 2, \dots , m-1$ and we get that  $\a$ is  $C^{2}$ in the variable $s,$  except when $s \in \{ \tau_1, \tau_2, \dots \tau_m  \},$ where $s=\tau_l(z, t), l = 1, 2, \dots, m$  is the point where $\a$ cuts the line $z=z_l$.

By the classical theory of smooth dependence of the solution of the ordinary differential equation on the data $(z,t)$, $\alpha(s,z,t)$ is $C^1$ in $s$ and the variables $z, t$, see Hartman \cite{Hartman} (Theorem 3.1, page 95). If $z_{m-1}<z<z_m,$ there is at most one point $s=\tau_m=\tau_m(z,t)$ where $\alpha(\tau,z,t)=0$. If no such $\tau$ exists then as before $\a(s,z,t)$ is $C^1$ in $s$ and in variables $z, t$ for all $0\leq s \leq t$. If such $\tau_m$ exists, because of the  hypothesis on $\l,$  $\alpha(s,z,t)$ is $C^1$ in $z, t$ for all $\tau_m \leq s \leq t$. Now the equation $\frac{ d \a}{ d s} = \l(\a, s)$ with initial data, $ \a(\tau_m(z, t), z, t)=z_m,$  has a unique solution for $ 0< t< \tau_m,$ which is $C^1$ in $s$ and which can be extended up to $t= \tau_m$  with right limit $\frac{ d \a( \tau_m +, z, t)}{ d s}$ and left limit $\frac{ d \a( \tau_m -, z, t)}{ d s}$.
Since $\l$ is $C^1$ in $z \geq z_m$, we have $\alpha(s,z,t) = \alpha(s,z_m, \tau_m(z, t))$ is  $C^1$ in $\tau_m(z, t)$ and $z_m$. Here again, we have used (Theorem 3.1, page 95 \cite{Hartman}). So if $\tau_m(z, t)$ is $C^1$ then $\a$ is $C^1$ in $z, t$.  We will now prove this.  \\

We first prove that $\tau_m=\tau_m(z, t)$ is $C^1$ in $(z, t)$ when $z <z_m. $ For $\tau_m \leq  s \leq t, \,\,\, \a(s, z, t)$  is strictly monotone in $s.$ So, it is invertible as a function of $s$ and it depends on the parameter $(z,t)$. Further, since $\a$ is $C^1$ in $(z, t),$ its inverse  ${\a}^{-1}$ is $C^1$ in $(z, t),$ by Inverse mapping theorem.  ${\a}^{-1}$ is $C^1$  from $[ \a(0,z,t), z_m]$ to $[0, \tau_m(z, t)]$, since it is extended to be $C^1$ in $\tau_m$. i.e ${\a}^{-1}$ is invertible and ${\a}^{-1}(z_m)=\tau_m(z, t).$ This is true up to $s= \tau_m(z, t)$  and hence $\tau_m(z, t)$ is $C^1$ in $(z, t)$.   In other words, since $\frac{ \partial \a}{ \partial \tau} \neq 0$  the equation $\alpha(s,z,t)=z_m$ has a unique solution $\tau_m$ which depends on  $(z, t)$ and 
$\frac{ \partial \tau_m}{ \partial z}= - \frac{ \partial \a}{ \partial z} \frac{1}{\l(0+, \tau(z, t))}, \,\,\, \frac{ \partial \tau}{ \partial t}= - \frac{ \partial \a}{ \partial t} \frac{1}{\l(0+, \tau_m(z, t))}.$  The same argument can be continued for the case when we solve characteristic equation \eqref{2.11} for $0<s<t $ with $(z,t)$ lying in the region $z_l<z<z_{l+1}, j=1,\dots,m-1$ and for the region $z<z_1$.  

The proof for the case $\lambda>0$ is similar.

{\bf  \flushleft Lemma 2.5:} Assume $\lambda(z,t)$ and $h(z,t)$ are $C^1$ in $z \neq z_l, l=1,\dots,m$ and their derivatives are bounded with their left and right limits existing at $z=z_l$. Suppose the initial data $v_0(z) \in C^{1}(\R)$. Then
there exists a unique $C^1$ solution $v(z, t)$ in $z \neq z_l$ for the problem
\begin{equation}
\begin{aligned}
&v_t + \l(z, t) v_z = h(z, t), \\
&v(z, 0)= v_0(z) 
\end{aligned}
\label{2.14}
\end{equation}
for $z \neq z_l$ with $v(z_l+, t) = v(z_l-, t)$ and it is given by 
\begin{equation}
v(z, t)= v_0( \a(0, z, t) ) + \int_0^{t} h(  \a( s, z, t), s) ds
\label{2.15}
\end{equation}
{ \bf Proof: }
Let $z \neq z_l,$ and let $\alpha(s,z,t)$ be the solution of
\[
\frac{d \alpha}{d s} = \l(\alpha, t), \,\,\, \alpha(s=t,z,t)= z,\,0<s<t.
\]
Along the curve $v=v(\alpha(s,z,t),s)$,
\begin{equation*}
\begin{aligned}
\frac{d v}{d s} &= \frac{\partial v}{\partial s}+ \frac{d \alpha}{d s}  \frac{\partial v}{\partial \alpha} = \frac{\partial v}{\partial s}+ \l(\alpha, t)  \frac{\partial v}{\partial \alpha}= h(\alpha, t)
\end{aligned}
\end{equation*}
where in the last equality, we used \eqref{2.14} with $\alpha=  \alpha(s,z,t), \,s$ is  the parameter. Integrating with respect to $s$, from $0$ to $t$ we get,  
\begin{equation*}
v(z, t)= v_0( \a(0, z, t) ) + \int_0^{t}{ h(  \a(s, z, t), s) } ds.
\end{equation*}
The smoothness property of $v(z,t)$  easily follows from the corresponding smoothness for $\alpha(s,z,t)$ given by the previous Lemma and differentiating the expression \eqref{2.15}. To do this we need to write the integral expression suitably. If $\alpha(s,z,t) \neq z_l$,  the characteristic curve does not have any intersection points with the lines $z=z_l$. Outside the lines $z=z_l$,  $\lambda(\alpha(s, z,t), t)$ and $\alpha(s, z,t)$ are smooth, then $v_j(z,t)$ is also smooth from the integral expression \eqref{2.15} and the regularity follows. If  for some time $s=\tau(z,t)$,  $\alpha(t,z,s) =z_l$, then write
\begin{equation*}
v(z, t)= v_0( \a(0, z, t) ) + \int_0^{\tau(z,t)}{ h(  \a(s, z, t), s) } ds +  \int_{\tau(z,t)}^{t}{ h(  \a(s, z, t), s) } ds
\end{equation*}
Now, we differentiate this expression with respect to $z$ and $t$ and use the smoothness property of $h$ for $z \neq z_l$ and Lemma 2.4 (ii) to get the result. 

To prove uniqueness, suppose there exists another solution $\tilde{v}.$ Then the difference $v_1= v - \tilde{v}.$ $v_1$ satisfies 
\begin{equation}
{v_1}_t + \l(z, t) {v_1}_z = 0, \,\,\, v_1(z, 0) =0
\label{2.16}
\end{equation} 
We have along the characteristic curve, $\frac{d v_1}{d t} = {v_1}_t + \l(z, t) {v_1}_z = 0$.  Integrating from $0$ to an arbitrary $t_0,$ we have,
$v_1(z_0, t_0) - v_1(\a(0, z_0, t_0), 0) = 0$ i.e. $v_1(z_0, t_0) =0.$ So, $v_1=0$ is the only solution to \eqref{2.16} and $v_1=0.$ Hence, $v = \tilde{v}$. \\

Now we complete the proof of the main theorem  (2.3) using Lemmas 2.4 (i) and (ii) and also Lemma 2.5. 

{\bf  \flushleft Proof of theorem 2.3.} As before for $z \neq z_l, $ let $A(z, t)$ be the matrix such that {
\[
B(z, t) = A(z, t)  \Lambda(z,t) A^{-1}(z, t)
\] 
where $\Lambda = diag( \l_1(z,t), \l_2(z,t), ...,\l_n(z,t))$. }In view of the assumptions $H_1$ and $H_2$, derivatives of $B(z,t), A(z,t)$ and $\lambda_j, j=1,\dots,n$ up to order two exist in $z\neq z_l, l=1,\dots,m$, and their right and left limits may not agree at $z=z_l$. 

Putting $u= A^{-1} v$, we transform \eqref{2.1} into 
\begin{equation}
v_t +\Lambda v_z = C v, \,\,\,v(z, 0) = v_0(z)
\la{2.17} 
\end{equation}
where $C = - A ( ( A^{-1})_t  + B ( A^{-1})_z ), \,\,\, \Lambda= diag( \l_1(z,t), \l_2(z,t), ...,\l_n(z,t))$  and $v_0(z)= Au_0(z)$.
Take $v = (v_1, v_2, v_3, \dots, v_n).$ We have 
\begin{equation}
\begin{aligned}
 \frac{\partial v_i}{\partial t}+ \l_i  \frac{\partial v_i}{\partial z}= \sum_{k=1}^n c_{ik} v_k
\end{aligned}
\la{2.18}
\end{equation}
where we used the matrix notation $C=(c_{ij})$. Using \eqref{2.15} in Lemma 2.5 and integrating \eqref{2.18} with respect to $s$ from $0$ to $t$ along the $j^{\rm th}$ characteristics we get, 
\begin{equation}
v_j(z, t) = v_{0j}( \a_j(0, z, t)) + \int_0^{t}{  \sum_{k=1}^n c_{jk}(\a_j(s, z, t), s) v_j(\a_j(s, z, t), s)     } ds
\label{2.19}
\end{equation}
We can write \eqref{2.19} as
\begin{equation*}
v = w_0 + Sv,
\end{equation*}
 where
vector $w_0$ has components $w_{0j}(z, t)=  v_{0j}( \a_j(0, z, t)) $, which is a fixed function, and
$S$ is the linear operator taking vector $v$ into a vector $w= Sv$ with components 
\begin{equation*}
w_j(z, t) = \int_0^{t}{  \sum_{k=1}^n c_{jk}(\a_j(s, z, t), s) v_j(\a_j(s, z, t), s)     } ds
\end{equation*}
For arbitrary fixed $T>0$, and  $L>0$, consider the space $C([-L,L]\times [0,T])$ of continuous  functions on 
$[-L,L]\times [0,T]$ with sup norm $||.||$. Since $\lambda_j(z,t)$ are bounded $\forall (z,t)$ there exists $\lambda_m<0$ and $\lambda_M>0$ such that $\lambda_m \leq \lambda_j(z,t) \leq \lambda_M$ $\forall j=1,\dots,n$ and $(z,t)$. So $\alpha_j(s,z,t) \in [-L+\lambda_m T,L+\lambda_M T]$, for all $0 \leq s \leq t$. Now, $S$ is a bounded linear map from $C([-L,L] \times [0,T])$ to itself with $|| S || = \sup_{|| v || =1}{ || S v ||}$,  bounded by $ T C_T$, where $C_T= \sup_{ [-L+\lambda_m T,L+\lambda_M T]\times [0,T], j,k} |c_{jk}((z, t)| $. So, the map $v \rightarrow w_0 +Sv$ is continuous and bounded. Further, it is a contraction map on $C([-L,L]\times [0,t_0]$ if $C_T .t_0<1$. By the fixed point theorem, the map
\[
v \rightarrow w_0 + Sv
\]
has a unique fixed point in $C([-L,L]\times [0,t_0],$ so that the integral equation \eqref{2.19} has a unique continuous solution for $0 \leq t_0$. This construction can be continued up to $0 \leq t \leq T$, since the characteristic curve $z=\alpha(s,z,t)$
exists $\forall (z,t)$ and $c_{jk}(z,t)$ is bounded by $C_T$ on $[-L+\lambda_n,L+\lambda_M] \times [0,T]$. Since $L$ and $T$ is arbitrary, there exists continuous $v=(v_1,v_2,...v_n)$ $\forall (z,t) \in R \times [0,\infty)$.

To prove that $v=(v_1,v_2...v_n)$ in \eqref{2.17}, is $C^1$ except on $z=z_j, j=1,\dots,m$, we write the components $v_j$  given by
the integral equation \eqref{2.16} as 
\begin{equation}
\begin{aligned}
v_j(z, t) &= v_j( \a_j(0, z, t))\\
& +{  \sum_{k=1}^n  \int_0^{\tau_{j1}(z,t)}c_{jk}(\a_j(s, z, t), s) v_j(\a_j(s, z, t), s)     } ds \\
&+\sum_{l=1}^{p_j} {  \sum_{k=1}^n \int_{\tau_{jl}(z,t)}^tc_{jk}(\a_j(s, z, t), s) v_j(\a_j(s, z, t), s)     } ds
\end{aligned}
\label{2.20}
\end{equation} 
where $\tau_{jl}(z,t), l=1,\dots,p_j$ is the time when the $j^{\rm th}$ characteristic starting at $(z,t)$ crosses the lines $z=z_l, j=1,\dots,m$ and at this time $\alpha_j(\tau_{jl}(z,t),z,t)=z_l$.
Differentiating \eqref{2.20} w.r.t. $z$ and $t$ and using $v$ is continuous and Lemma 2.5, we get $v$ is continuously differentiable in $z\neq z_l$. Hence, $v$ is continuous in $R \times [0,\infty)$ and satisfies the equation and initial condition in \eqref{2.16} for $z \neq z_l$. Then, $u=A v$ solves the problem $\eqref{2.1}$ with interface condition \eqref{2.4}, i.e. $A^{-1} u$ is continuous across $z=z_l$. 

Uniqueness follows as in the proof of Lemma 2.4, by taking a difference between two solutions and showing the solution of the homogeneous equation with zero function as the initial data has only a trivial solution. 

\section{ Symmetric hyperbolic system and interface conditions  } 
In this section, we study the initial value problem for  a symmetric hyperbolic system of the form  
\begin{equation}
Lu=B_0(x) u_t + \sum_{j=1}^n B_j(x) u_{x_j} = f(x,t),
\label{3.1}
\end{equation}
in the region $R_T=R^n \times [0,T],$ with initial conditions
\begin{equation}
u(x,0)=u_0(x).
\label{3.2}
\end{equation}The unknown $u$ is a mapping from $R_T$ to $R^m$. We assume the $m\times m$ matrices $B_j$ are at least once continuously differentiable except along the
surface $\{(x_1,\dots,x_n) : x_n=0\}$. Since  we need to consider the half spaces  $\{x = { \it (x_1,\dots,x_n)  :x_n > }0 \}$ and  $\{x= { \it (x_1,\dots,x_n)} :x_n<0 \}$,  it is convenient to single out the variable $x_n$ and write the 
space variable $x=(x',x_n)$ where $x'=(x_1,,,x_{n-1})$ and 
$\tilde{R}_T={R_T}\backslash {\Gamma}_T,$ where ${\Gamma}_T= \{  x=(x',x_n) : x_n =0  \} \times [0, T).$ We require the equation to be satisfied in  $\tilde{R}_T$ with Interface condition on ${\Gamma}_T,$ satisfied in the weak sense that we shall see in Theorem 3.3 here.   \\ 

{\bf  \flushleft Assumptions on $B_j(x), j=0,1,...n$:}\\
A1 : There exists positive constants $c_1$ and $c_2$ such that $c_1 I \leq B_0 \leq c_2 I$ , where $I$ is the identity matrix of order $m$. \\
A2 : $B_j(x), j=0,1,...n,$ are symmetric $n \times n$ matrices whose entries are once continuously differentiable except along the hypersurface ${\Gamma}_T$ and the traces  $B_j(x',0\pm)$ exists.\\
A3 : There exists $C>0, $ such that for $j=0,\dots,n$
\[
\begin{aligned}
 &||(B_{j})||_{L^\infty(R^n)}  \leq C,j=1,\dots,n,\,\,       ||(B_{j})_{x_j}||_{L^{\infty}(R^n)}  \leq C,  j=1,\dots,n-1\\
 &|| (B_n)_{x_n})||_{L^{\infty}(R^{n-1}\times (0, \infty)} \leq C, \,\, ||(B_{n})_{x_n}||_{L^{\infty}(R^{n-1}\times(- \infty, 0)}  \leq C,
\end{aligned}
\] 

We are interested in the weak solution of the initial value problem for \eqref{3.1} with initial conditions \eqref{3.2}
in a suitable function space. See \cite{Evans}. One of the most important steps in the existence of solutions in the $L^2$ Sobolev space is a derivation of the energy estimates. Here, since the coefficient matrices are not smooth along the hyperplane  $x_n=0$, we need an interface condition along the surface of discontinuity to get the energy estimate and the weak formulation of the solution globally. To derive energy estimates it is more convenient to work with the following norm
\[
||| u(t) ||| = (\int_{R^n}{( B_0(x) u(t), u(t)) dx })^{1/2}
\] 
 Observe that because of the assumption $A_1$, for $u=(u_1,...u_m)$ an unknown function depending on $(x,t)$, the map $ ||| \cdot ||| \rightarrow [0,\infty)$ defined above 
is a  norm equivalent to the norm $|| u(t) ||_{L^2} = (\int_{R^n}{(u(t), u(t))dx })^{1/2}$ as a function of $x,$ $\forall$ fixed $t$. 

In Section 3.1 we get the weak formulation for the initial value problem to system \eqref{3.1} with initial condition \eqref{3.2} and derive the interface condition on $\Gamma_T.$ The energy estimate that we use to get the existence result for our system, is derived in Section 3.2. We then prove the existence of a unique solution to our system in Section 3.3.

\subsection{ Weak formulation for symmetric hyperbolic system and interface conditions  }

To formulate the weak solution we introduce the following function spaces. {
\begin{equation}
\begin{aligned}
H=\{ v \in L^2(R^n) :& v|_{R^{n-1}\times (-\infty,0)} \in H^1(R^{n-1}\times (-\infty,0)),\\
&v|_{R^{n-1}\times (0,\infty)} \in H^1(R^{n-1} \times (0,\infty))\}
\end{aligned}
\label{3.3}
\end{equation}
}
For each element $v$ in $H$ the traces $v(x',0+), v(x',0-)$ exists. We define 
\begin{equation}
\hat{H}=\{ v \in H : B_n(x',0+)v(x',0+)=B_n(x',0-)v(x',0-) \}
\label{3.4}
\end{equation}
For $u \in H,$ multiply $Lu$ in \eqref{3.1} by  $\phi$, with $\phi(x,T)=0,  \phi \in C^1(R_T;R^m)$  and integrate with respect to $(x,t)$ in the region $R_T$, we get
\begin{equation}
\int_0^T\int_{R^n} (Lu, \phi )dx dt = \int_0^T\int_{R^n} (B_0(x) u_t+ \sum_{j=1}^n B_j(x) u_{x_j} , \phi) dx dt
\label{3.5}
\end{equation}
Considering each of the terms in this expression, using the symmetry of $B_j, j=0,1..,n$, and integrating by parts we get the following equalities
\begin{equation}
\begin{aligned}
 \int_0^T\int_{R^n} (B_0(x) u_t, \phi) dx dt =& - \int_0^T\int_{ R^n} (u, (B_0(x)\phi)_t) dx dt \\
&-\int_{R^n} (B_0(x) u(x,0), \phi(x,0)) dx,
\end{aligned}
\label{3.6}
\end{equation}
\begin{equation}
\begin{aligned}
 \int_0^T\int_{R^n} (B_j(x) u_{x_j}, \phi )dx dt &=\int_0^T\int_{R^n} (u_{x_j},B_j \phi )dx dt \\
&= - \int_0^T\int_{R^ n} (u, (B_j(x)\phi)_{x_j} )dx dt , \\
\end{aligned}
\label{3.7}
\end{equation}
{ for $j=1,\dots,n-1$}, and 
\begin{equation}
\begin{aligned}
 \int_0^T&\int_{R^1} (B_n(x) u_{x_n}, \phi )dx dt = \int_0^T\int_{R^1} (u_{x_n}, B_n(x)\phi )dx_n dt\\
&=\int_0^T \int_{-\infty }^0 (u_{x_n},B_n \phi )dx dt+\int_0^T\int_0^{\infty} (u_{x_n},B_n \phi )dx_n dt \\
&=(B_n(x',0-) u(x',0-,t),\phi(x',0,t))-\int_{\infty}^0(u,(B_n(x)\phi)_{x_n})dx_n\\
&+(B_n(x',0+) u(x',0+,t),\phi(x,0,t))-\int_0^{\infty} (u,(B_n(x)\phi)_{x_n})dx_n\\
&=(B_n(x',0+) u(x',0+,t)-B_n(x',0-) u(x',0-,t),\phi(x,0,t))\\
&-\int_{R^1} (u,(B_n(x)\phi)_{x_n})dx_n
\end{aligned}
\label{3.8}
\end{equation}
In \eqref{3.7}, \eqref{3.8} the integral terms on the boundary of the region vanish as  $\phi(x, T)=0$. Using \eqref{3.6}-\eqref{3.8} in \eqref{3.5}, we get
\begin{equation}
\begin{aligned}
\int_0^T&\int_{R^n}(Lu, \phi )dx dt= \\
&-\int_{R^n} (B_0(x) u(x,0), \phi(x,0)) dx - \int_0^T\int_{R^n } (u, L^T\phi) dxdt \\
&-\int_0^T \int_{R^{n-1}}[(B_n(x',0+) u(x',0+,t)\\
&-(B_n(x',0-) u(x',0-,t),\phi(x',0,t))]dx' dt
\end{aligned}
\label{3.9}
\end{equation}
where $L^T$  is  defined by
\begin{equation}
L^T u= -[(B_0(z)u)_t +\sum_{j=1}^n (B_j(x)u)_{x_j}]
\label{3.10}
\end{equation}
and the interface condition is 
\begin{equation}
B_n(x',0+) u(x',0+,t)-(B_n(x',0-) u(x',0-,t) =0.
\label{3.11}
\end{equation}
Moreover, for $u \in \hat{H}$  in \eqref{3.4} and $\phi \in  C_c^1(R^n \times (0,T))$, we have 
\begin{equation*}
\int_0^T \int_{R^n}(Lu,\phi) dx dz dt= \int_0^T \int_{R^n}  (u,L^T \phi),\,\, \,\,\, \forall \phi \in  C_c^1(R^n \times (0,T))
\end{equation*}
and $L^T$ is called the transpose of $L$. This is using \eqref{3.11} and \eqref{3.9}. Now using the above equation, we introduce the notion of weak solution to system \eqref{3.1}, \eqref{3.2}.  \\

{\bf  \flushleft Definition :}
We say $u\in C([0,T],L^2(R^n))$ is a weak solution of \eqref{3.1}, \eqref{3.2} if for all $\phi \in  C^1(R_T,R^m)$ with $\phi(x,T)=0,$ for $T>0$ 
\begin{equation*}
\begin{aligned}
\int_0^T \int_{R^n} (f,\phi) dxdt+\int_{R^{n}}(B_0(x) u_0(x), \phi(x,0)) dx=
\int_0^T \int_{R^n}(u,L^T \phi)dx
\end{aligned}
\end{equation*}
To prove that there exists a unique solution to the Cauchy problem to this system we need Energy estimates that we derive in the following section. 

\subsection{Energy estimates for symmetric hyperbolic system and interface conditions} 

Here we get the energy estimates for our system \eqref{3.1}, \eqref{3.2}. We introduce a norm that is equivalent to the $L^2$ norm, which is the same norm as defined at the beginning of Section 3, by: 
\begin{equation}
|||u|||= \Big(\int_{R^n}\displaystyle{(B_0 u,u) dx\Big )}^{1/2}%{\displaystyle{\frac{1}{2}}}
\label{3.14}
\end{equation}
The fact that the norm in \eqref{3.14} is equivalent to $L^2$ norm, follows from our assumption  $A1$.

{\bf  \flushleft Lemma 3.1} Under the assumptions A1, A2, A3, $\forall u\in C([0,T],{H})$
satisfying 
\begin{equation*}
Lu=f
\end{equation*}

\eqref{3.1} in $z\neq 0$ we have the estimate
 \begin{equation*}
 \begin{aligned}
& \frac{d}{dt}( e^{- C(n) t} ||| u(t) |||^2 ) \leq e^{-C(n)t}||f(t)||_{L^2(R^n)}^2\\
&+\frac{1}{2} e^{- C(n) t}   \int_{R^{n-1} }[(B_{n}(x',0+) u(x', 0+, t), u(x', 0+, t)) \\
&- (B_{n}(x',0-) u(x', 0-, t), u(x', 0-, t))]dx' 
\end{aligned}
\end{equation*}
where $C(n)=1+nC$ with $C$ as in assumption $A3$,  \\
which is obtained using the equalities
\begin{equation*}
\int_{R^ n} (B_j(x)u,  u_{x_j} )dx =- \frac{1}{2}\int_{R^ n} (u, (B_j(x)_{x_j}u )dx 
\end{equation*}
for $j=1,\dots,n-1$ and 
\begin{equation*}
\begin{aligned}
 \int_{R^1}& (B_n(x) u_{x_n}, u )dx =- \frac{1}{2}\int_{R^1} (u, (B_n(x))_{x_n}u)dx_n \\
&=-\frac{1}{2}(B_n(x',0+) u(x',0+,t),u(x',0+,t))\\
&+\frac{1}{2}(B_n(x',0-) u(x',0-,t),u(x',0-,t))
\end{aligned}
\end{equation*}
 for $j=n$.
 \\
{\bf Proof : }  Since we consider $u\in C([0,T],{H})$, we write $u(x,t)$ as $u(t)$, a function of $t$ valued in $H$ defined in \eqref{3.3}  and may also suppress the variable $(x,t)$ and write $u$ in our analysis. We find
\begin{equation}
\begin{aligned}
\frac{d}{dt}( \frac{1}{2} {|||u|||}^{2} ) &= \frac{1}{2} \frac{d}{dt}
\int_{R^n}(B_0(x) u, u)dx   \\
&=  \int_{R^n }(B_0(z) u_t, u)dx \\
&= -\sum_{j=1}^n \int_{R^n} (B_j(z) u_{x_j}, u) dx +\int_{R^n} (f,u)dx.
\end{aligned}
\label{3.15}
\end{equation}
Now for $j=1,\dots,n-1$, integrating by parts w.r.t. $x$ and using symmetry of $B_j$, we get 
\begin{equation*}
\begin{aligned}
 \int_{R^n} &(B_j(x) u_{x_j}, u )dx  =\int_{R^n} (u_{x_j},B_j u)dx  \\
&= - \int_{R^ n} (u, (B_j(x) u)_{x_j} )dx\\
&=- \int_{R^ n} (u, B_j(x) u_{x_j} )dx  - \int_{R^ n} (u, (B_j(x)_{x_j}u )dx  \\
&=- \int_{R^ n} (B_j(x)u,  u_{x_j} )dx  - \int_{R^ n} (u, (B_j(x)_{x_j}u )dx  
\end{aligned}
\end{equation*}
Taking the first term in the RHS of the above equation to left we get,  
\begin{equation}
\int_{R^ n} (B_j(x)u,  u_{x_j} )dx =- \frac{1}{2}\int_{R^ n} (u, (B_j(x)_{x_j}u )dx 
\label{3.15}
\end{equation}

For $j=n$, we get
\begin{equation}
\begin{aligned}
 \int_{R^1} &(B_n(x) u_{x_n}, u )dx_n = \int_{R^1} (u_{x_n}, B_n(x)u)dx_n\\
&= \int_{-\infty }^0 (u_{x_n},B_n u)dx_n +\int_0^{\infty} (u_{x_n},B_n u)dx_n  \\
&=(B_n(x',0-) u(x',0-,t),u(x',0-,t))-\int_{\infty}^0(u,(B_n(x)u)_{x_n})dx_n\\
&-(B_n(x',0+) u(x',0+,t),u(x',0+,t))-\int_0^{\infty} (u,(B_n(x)u(x)_{x_n})dx_n  \\
&=-\int_{R^1} (u,(B_n(x)u)_{x_n})dx_n-(B_n(x',0+) u(x',0+,t),u(x',0+,t))\\
&+(B_n(x',0-) u(x',0-,t),u(x',0-,t))
\end{aligned}
\label{3.16}
\end{equation}

Here again, carrying out the differentiation in the integrand in the last term in \eqref{3.16} and using the symmetry of $B_n$, we get as before.

\begin{equation}
\begin{aligned}
 \int_{R^1}& (B_n(x) u_{x_n}, u )dx_n =- \frac{1}{2}\int_{R^1} (u, (B_n(x))_{x_n}u)dx_n \\
&=-\frac{1}{2}(B_n(x',0+) u(x',0+,t),u(x',0+,t))\\
&+\frac{1}{2}(B_n(x',0-) u(x',0-,t),u(x',0-,t))\\
\end{aligned}
\label{3.17}
\end{equation}
Integrating \eqref{3.17} with respect to $(x_1,,,x_{n-1})$  we get 
\begin{equation}
\begin{aligned}
 \int_{R^n}(B_n(x) u_{x_n}, u )dx &= -\frac{1}{2}\int_{R^n} (u, (B_n(x))_{x_n}u)dx\\
&=-\frac{1}{2}\int_{R^{n-1}}(B_n(x',0+) u(x',0+,t),u(x',0+,t))dx'\\
&+\frac{1}{2}\int_{R^{n-1}}(B_n(x',0-) u(x',0-,t),u(x',0-,t)dx' )\\
\end{aligned}
\label{3.18}
\end{equation}
From \eqref{3.14}-\eqref{3.18}, we have
\begin{equation*}
\begin{aligned}
\frac{d}{dt}( \frac{1}{2} {||| u(t) |||}^{2} ) &= \sum_{j=1}^n\frac{1}{2}\int_{R^1} (u, (B_j(x))_{x_j}u)dx+
\int_{R^n} (f,u)dx\\
&+\frac{1}{2}\int_0^T\int_{R^{n-1}}(B_n(x',0+) u(x',0+,t),u(x',0+,t))dx'dt\\
&-\frac{1}{2}\int_0^T\int_{R^{n-1}}(B_n(x',0-) u(x',0-,t),u(x',0-,t)dx' dt\\
\end{aligned}
\end{equation*}
By assumption $A3$ and Cauchy Schwarz inequality, we have 
\begin{equation}
\begin{aligned}
 \frac{d}{dt}( \frac{1}{2} {||| u(t) |||}^{2} ) &\leq n C ||u||_{L^2(R^n)}^2 + ||f(t)||_{L^2(R^n)}^2 +  ||u||_{L^2(R^n)}^2 \\
&+\frac{1}{2}\int_{R^{n-1}}(B_n(x',0+) u(x',0+,t),u(x',0+,t))dx'\\
&-\frac{1}{2}\int_{R^{n-1}}(B_n(x',0-) u(x',0-,t),u(x',0-,t)dx' \\
\end{aligned}
\label{3.20}
\end{equation}
This inequality \eqref{3.20} can be written as 
\begin{equation}
\begin{aligned}
 \frac{d}{dt}&( \frac{1}{2} e^{-C(n)t}{||| u(t) |||}^{2} ) \leq e^{-C(n)t} ||f(t)||_{L^2(R^n)}^2 \\
&+e^{-(C(n))t}[\frac{1}{2}\int_{R^{n-1}}(B_n(x',0+) u(x',0+,t),u(x',0+,t))dx'\\
&+\frac{1}{2}\int_{R^{n-1}}(B_n(x',0-) u(x',0-,t),u(x',0-,t)dx'] \\
\end{aligned}
\label{3.22}
\end{equation}
From this Lemma, the following result is immediate using \eqref{3.22}. 
{\bf  \flushleft Corollary 3.2 :} $\forall u \in C^1([0,T],H)$ with $u(x,T)=0$, and 
\begin{equation*}
\begin{aligned}
\int_{R^{n-1}}[(B_n(x',0+) &u(x',0+,t),u(x',0+,t))\\
&-(B_n(x',0-) u(x',0-,t),u(x',0-,t))]dx' \leq 0,
\end{aligned}
\end{equation*}
we have
\begin{equation}
\begin{aligned}
 ||| u(t) |||^{2}  &\leq C_T ||Lu||_{L^2(R^n \times [0,T])}^2 \\
\end{aligned}
\label{3.24}
\end{equation}

\subsection {Existence result for the initial value problem} 
In this section, we prove the existence of the Cauchy problem in the space-time region $R_T=R^n \times [0, T]$.
Note that for the weak formulation   of our solution we need the interface conditions 
\begin{equation}
B_n(x',0+)v(x',0+)=B_n(x',0-)v(x',0-)
\label{3.25}
\end{equation}
and for energy inequality, we need
\begin{equation}
\begin{aligned}
\int_{R^{n-1}}&[(B_n(x',0+) u(x',0+,t),u(x',0+,t))\\
&-(B_n(x',0-) u(x',0-,t),u(x',0-,t)]dx' \leq 0 .
\end{aligned}
\label{3.26}
\end{equation}
So for an existence theory based on these,  the solution space must satisfy both these conditions  \eqref{3.25} and \eqref{3.26}  simultaneously. A sufficient condition for this requirement is\\
\begin{equation}
\text{either }\,\,u(x', 0\pm)=0\,\, \text{or} \,\,B_n(x', 0+)\leq 0,\,\, \text{and} \,\,B_n(x',0-) \geq 0.
\label{3.27}
\end{equation}
 Depending upon the conditions required for energy estimates we define different solution spaces.   When $B_n(x',0+) \leq 0$ and $B_n(x',0) \geq 0$ define the function space
\[
\begin{aligned}
\tilde{C^1}(R_T)=\{&v \in C^1(R_T,R^m) : v(T)=0,\\
& B_n(x',0+)v(x',0+)=B_n(x',0-)v(x',0-)\\
&supp (v(.,t)),\,\, compact\,\, in \,\,R^n , 0\leq t\leq T\}
\end{aligned}
\]
For general $B_n$,  define
\[
\begin{aligned}
\tilde{C_1^1}(R_T)=\{&v \in C^1(R_T,R^m) : v(T)=0, v(x',0\pm,t)=0,\\
& supp (v(.,t)),\,\, compact\,\, in \,\,R^n , 0\leq t\leq T\}
\end{aligned}
\]

Define an inner product in $\tilde{C^1(R_T)}$
\[
(v,u)_{\cal{H}}=\int_0^T \int_{R^n}(L^Tv,L^Tu)dxdt
\]
Let $\cal{H}$ be the completion of $\tilde{C^1(R_T)}$ with respect to the norm induced by this inner product   
\[
||v||_{\cal{H}}=(\int_0^T\int_{R^n} |L^Tu|^2 dxdt)^{1/2}.
\]
 where $L^Tu$ defined in \eqref{3.10} is adjoint operator of $Lu$.

{\bf  \flushleft Theorem 3.3} Under the assumptions $A1,A2,A3$ we have the following.
\begin{itemize}
\item Given $f \in L^2(R_T)$ and $u_0 \in L^2(R^n)$ there exists a unique weak solution $u \in C([0,T], L^2(R^n))$ of \eqref{3.1} and \eqref{3.2} satisfying the interface condition 
\[
\begin{aligned}
({L^T})^{-1}u(x',0) &\in   \text{Closure}\{v \in C^1(R_T,R^m) : v(T)=0, v(0\pm,t)=0,\\
& supp (v(.,t)),\,\, compact\,\, in \,\,R^n , 0\leq t\leq T\}
\end{aligned}
\]
\item If $B_n(x',0-) \leq 0$ and $B_n(x',0-)\geq 0$ then there exists a unique weak solution $u \in C([0,T], L^2(R^n))$ of \eqref{3.1} and \eqref{3.2} satisfying the interface condition 
\[
\begin{aligned}
({L^T})^{-1}u(x',0) &\in \text{Closure} [\{v \in C^1(R_T,R^m) : v(T)=0,\\
& B_n(x',0+)v(x',0+)=B_n(x',0-)v(x',0-)\\
&supp (v(.,t)),\,\, compact\,\, in \,\,R^n , 0\leq t\leq T\}]
\end{aligned}
\]
\end{itemize}
{\bf Proof:} Recall that the energy estimate \eqref{3.24} shows that, $\forall v\in H$
\begin{equation*}
(\int_0^T\int_{R^n} |v|^2 dx dt)^{1/2} \leq C_T (\int_0^T\int_{R^n} |L^Tv|^2 dxdt)^{1/2}.
\end{equation*}
Now consider the linear functional $F: {\cal{H}} \rightarrow R^1$
\[
F(v)=\int_{R^n} (u_0,v(x,0))dx +\int_0^T\int_{R^n}(f,v)dxdt
\]
By Cauchy Schwarz inequality, we have 
\[
\begin{aligned}
|F(v)| & \leq ||u_0||_L^2(R^n) ||v(x.0)||_{L^2(R^n)} + ||f||_{L^2(R_T)}  ||v||_{L^2(R_T)}\\
& \leq ( ||u_0||_L^2(R^n) + ||f||_{L^2(R_T)})  ||L^Tv||_{L^2(R_T)}
\end{aligned}
\]
So $F$ is a continuous linear function on ${\cal{H}}$. By Riesz representation theorem there exists a unique $u_1 \in {\cal{H}}$
such that
\[
F(v)=(u,v)_{\cal{H}}
\]
This means that there exists $u_1 \in H$ such that $\forall v \in {\cal{H}}$
\[
\int_{R^n} (u_0,v(x,0))dx +\int_0^T\int_{R^n}(f,v) dx dt =\int_0^T \int_{R^n}(L^Tu_1,L^Tv) dx dt
\]
Now take $u=L^Tu_1$, clearly $u \in L^2(R_T)$ and satisfies
\begin{equation}
\int_0^T \int_{R^n}(L^Tu_1,L^Tv)dxdt=\int_{R^n} (u_0,v(x,0))dx +\int_0^T\int_{R^n}(f,v)dxdt 
\label{3.28}
\end{equation}
$\forall v \in {\cal{H}}$

So \eqref{3.28} proves that $u$ is a weak solution in $L^2$.

Next, we come to the interface condition at $x_n=0$. By the energy estimate \eqref{3.24},  $L^T$ is a continuous linear functional from $H$ to $L^2(R_T)$ and is invertible. It follows that $u$ satisfies the interface conditions.

{\bf  \flushleft Remark:}  In this Section, we are interested in the existence of a solution with interface conditions. 
Here we make a few remarks.
\begin{itemize}
\item We have the existence of a unique solution for the 1D case in Section 2 ( with no requirement for the symmetry condition) using the method of characteristics with interface conditions \eqref{2.4}. We have the existence of a unique solution for the nD case for symmetric hyperbolic systems with interface conditions \eqref{3.27}. An interesting question is whether there is a comparison between the two interface conditions obtained by the method of characteristics and by the energy method to prove the existence theorem in section 3.

\item The solution we constructed by the method of characteristics satisfies an interface condition that is not the same, in general as the interface condition by energy method that we used in section 3 to prove the existence of a solution. Comparing these two interface conditions is not feasible as we see below.

Here is the method of characteristics argument, waves from $z<0$ and $z>0$ cross over and interact. On the other hand, the above existence theorem in Section 3 does not allow wave crossing. So it is important to analyze the problem
by a suitable approximation method, which we have done in this paper.  \\

\end{itemize}

{\section{Conclusion}
\la{Sum-sect} 

In this paper, we are analyzing wave propagation in a piece-wise homogeneous medium. Mathematically, this amounts to studying the initial value problem for a system of linear hyperbolic equations with a coefficient matrix having a discontinuity along a surface. For a strictly hyperbolic system in one space variable, we use the method of characteristics to construct the solution. The solution constructed this way is discontinuous across the surface of discontinuity of the coefficient matrix but the components in the eigen-directions are continuous. Here we remark that the characteristic speed and the corresponding eigenvectors are discontinuous along the surface. Next, we consider another method of construction of solutions using energy methods. We consider a symmetric hyperbolic system and derived the interface condition to get energy estimates. For the one-dimensional case with the additional condition of symmetry, we compared the interface condition. We get that the interface condition obtained by the method of characteristics satisfies the one obtained by the energy method if the characteristic speed on the right of the discontinuity is smaller than that of the left. 

The significance of this paper is that in the previous work in \cite{Fatemeh} the matrices involved were symmetric and one constant for the region $-1 < x < 0, t> 0$ and another constant matrix of a possibly different order, in the region $0 < x < 1, t> 0$  but in our work the matrix involved $B$ is a function of $z$ having a discontinuity on $z=0.$ The focus in \cite{Fatemeh} is the numerical analysis for their system whereas we have a more theoretical approach to our analysis. We do consider the case when the matrix $B$ is piecewise constant with a discontinuity on $z=0$ in Section 2.1  and write explicitly the solution to this system with given initial conditions. We also do the same for  $B$ having  2 discontinuities along $z=z_1, \,\,\, z=z_1$.  This analysis for when  $B$ having $m$ discontinuities along $m$ curves  $z=z_l, \,\,\, l = 1,\dots, m,$ $m$ is any natural number is the same as the analysis done here for $m=2$. The general system in 1D - space variable in the domain $\{(z,t): z \in \R^1, t>0\}$ with $B$ having discontinuities along $m$ curves  $z=z_l, \,\,\, l = 1,\dots, m$  is studied in Section 2.2 and we prove the existence and uniqueness of the solution to this corresponding system. In our work, the matrix $B$ need not be symmetric. However, symmetricity is required in Section 3, where we derive energy estimates that we use to construct solutions to the initial value problem for the $n$ - space dimensional symmetric hyperbolic system with discontinuity along the surface $\{(x_1,x_2, \dots,x_n): x_n=0\}$ and derive the corresponding interface condition. Such an analysis has not been done previously as per our knowledge.  

As a continuation of our work, we can follow the study by C. L. Pekeris \cite{Pekeris}.} It is possible to do the same analysis {that is done in this work of ours} for $[0,\infty)$ as well as for $[0,1]$. We expect similar results would hold. Technically, one needs to take into account the boundary conditions in the integral equations in Lemma 2.5.

{\section{Acknowledgements}
\la{Sum-sect} 

I would like to thank Professor Boris Belinskiy, University of Tennessee at Chattanooga, Tennesse, USA  for suggesting this problem and for his valuable comments, support that helped complete this paper. \\

{ \flushleft Data availability statement: None \\ 
Funding statement: This research did not receive any specific grant from funding agencies in the public, commercial, or not-for-profit sectors. \\
Conflict of interest disclosure: None \\
}

\section*{ORCID}
Kayyunnapara Divya Joseph https://orcid.org/0000-0002-4126-7882


\begin{thebibliography}{99}


\bibitem{Boris}
    \newblock  Belinsky  B. P, Andronov  I.V.
    \newblock   On Boundary -- Contact Acoustic Problems for a Vertically Stratified
Medium Bounded from Above by a Plate with Concentrated Inhomogeneities.
    \newblock PMM USSR (J. of Applied Mathem. and Mechanics).  { \bf 54} no. 3, (1990) 366-371.


\bibitem{Christodoulou}
\newblock  Christodoulou D, O'Murchadha  N.
\newblock The Boost Problem in General Relativity.
\newblock Commun. Math. Phys.  { \bf 80} (1981) 271-300.


\bibitem{Courant}
    \newblock  Courant R,  Hilbert D.
\newblock Method of Mathematical Physics.
  \newblock  Vol. II (1989) Chapter 5, Section 6.2.



\bibitem {Evans} 
\newblock Evans  L. C. 
\newblock Partial differential equations.
\newblock AMS. (1998) 

\bibitem{Dym}
    \newblock  Dym H.
\newblock Linear Algebra  in Action. 
\newblock   AMS. (2006) Vol. 78, Chapter 4, Section 4.6.

\bibitem{Fatemeh}
    \newblock Ghasemi F, Nordstrom J. 
    \newblock Coupling Requirements for Multiphysics Problem posed on two domains.
    \newblock  SIAM J. Numer. Anal.   Vol. 55, no. 6 (2017) 2885-2904.                                                                                                                                                                                                                                                                                                                                                                                                                                                                                                                                                                                                                                                                                                                                                                                                                                                                                                                                            

\bibitem{Lax3}
\newblock  Friedrichs K.O, Lax  P. D.,
\newblock  Systems of conservation equations with a convex extension.
\newblock Proc. Natl. Acad. Sci. USA. { \bf 68} (1971) 1686-1688.

\bibitem{Fredrichs1}   
\newblock Friedrichs K.O. 
\newblock Symmetric hyperbolic linear differential equations.
\newblock Comm. Pure Appl. Math. { \bf 7} (1954)  345-392.

\bibitem{Friedrichs2}
\newblock  Friedrichs K.O.
\newblock Symmetric positive linear differential equations.
\newblock  Comm. Pure Appl. Math. Vol. 11  (1958) 333-418.

\bibitem{Friedrichs3}
\newblock Friedrichs  K.O.
\newblock Nonlinear Hyperbolic Differential Equations for Functions of Two Independent Variables
\newblock  Am. Jl. Math.   Vol. 70, no. 3 (1948) 555-589.

\bibitem {g1}
\newblock Garding L. 
\newblock Linear hyperbolic partial differential equations with constant coefficients.
\newblock Acta Math.{ \bf 85} (1951) 1-62. 

\bibitem {h1}  
\newblock Hadamard  J. 
\newblock Lectures on Cauchy Problem in Linear Partial Differential equations.
\newblock Yale University press. (1923).

\bibitem{Hartman}
    \newblock Hartman P. 
\newblock Ordinary Differential Equations.
\newblock SIAM, Second Edition.

\bibitem{Hoffman}
\newblock  Hoffman K, Kunze  R. 
\newblock Linear Algebra.
\newblock PHI Learning Private Limited. Second Edition, Theorem 2, 187.


\bibitem{Hughes}
\newblock Thomas J . R. Hughes, Tosio Kato, Jerrold E. Marsden.
\newblock Well-posed Quasi-linear Second-order Hyperbolic Systems with Applications to Nonlinear Elastodynamics and General Relativity.
\newblock Arch. Rat. Mech. Anal. no.3, { \bf 63} (1977) 273-294.

\bibitem{John}
    \newblock  John F.
\newblock Partial differential equations.
\newblock Springer. 4th Edition, Chapter 2, Section 5.

\bibitem{KTJ}
\newblock Joseph  K. T. 
\newblock A Riemann problem whose viscosity solution containing $\delta$-measures.
\newblock Asym. Anal. { \bf 7} (1993) 105-120.


\bibitem{Kopriva}
\newblock Kopriva D. 
\newblock A spectral multidomain method for the solution of hyperbolic systems.
\newblock NASA Contractor Report 178105. ICASE. (1986).




\bibitem{Kreiss}
\newblock  Kreiss H. O.
\newblock  Problems with different time scales for partial differential equations.
\newblock Comm. Pure Appl. Math. { \bf 33} (1980) 399-439.


\bibitem{Lax}
    \newblock  Lax P. D.
\newblock Linear Algebra and its Applications.
\newblock John Wiley and Sons. Second Edition. (2007) 130.

\bibitem{Lax2}
\newblock  Lax P. D,  Phillips R. S. 
\newblock Local boundary conditions for dissipative symmetric linear differential operators.
\newblock Comm. Pure Appl. Math.,  Vol. 13, Issue 3 (1960) 427-455.


\bibitem{Lefloch}
\newblock  LeFloch P. G.
\newblock An existence and uniqueness result for two nonstrictly  hyperbolic systems, in  Nonlinear Evolution Equations that change  type, (eds) Barbara Lee Keyfitz and Michael Shearer.
\newblock  IMA, { \bf 27} (1990) 126-138.

\bibitem{Leray}
\newblock  Leray J.
\newblock Lectures on hyperbolic partial differential equations.
\newblock  Institute for Advanced Study. (1953).




\bibitem {Ya} 
\newblock Ya. B. Lopatinskii.
\newblock The mixed Cauchy type problem for the equation of hyperbolic type.
\newblock Dopovidi Akad. Nauk. Ukrain. RSR Ser. A, (1970) 592-594.


\bibitem{Pekeris}
\newblock  Pekeris C. L.
\newblock Theory of Propagation of Sound in a Half-Space of Variable Sound Velocity under Conditions of Formation of a Shadow Zone.
\newblock The Journal of the Acoustical Society of America.{ \bf 18} (1946) 295 



\bibitem {p1} 
\newblock Petrowsky  I.
\newblock Uber das Cauchysche Problem fur Systeme von partiellen Differentialgleichungen.
\newblock Rec. Math. [Mat. Sbornik].  Vol. 2, no. 5, { \bf 44} (1937) 815-870.



\bibitem{Quarteroni}
\newblock Quarteroni A, Lassila  T, Rossi S, Ruiz-Baiera R.
\newblock Integrated Heart - Coupling multiscale and multiphysics models for the simulation of
the cardiac function.
\newblock Computer Methods in Applied Mechanics and Engineering.  { \bf 314} (2017) 345-407.


\bibitem{Rauch}
\newblock   Rauch J. B, Massey III F. J.
\newblock Differentiability of Solutions to Hyperbolic Initial Boundary Value Problems.
\newblock Transactions of the American Mathematical Society.  Vol. 189 (1974) 303-318.

\bibitem {s1} 
\newblock Schauder  J.
\newblock Das Anfangswert problem einer quasilinearen hyperbolischen  Differentialgleichung zweiter Ordnung in beliebiger Anzahl von un abhahangigem.
\newblock Ver-nderlichen. Fund. Math. { \bf 24 } (1935) 213-246.


 \end{thebibliography}
\end{document}